\documentclass[10pt]{article}
\usepackage{cite}
\usepackage{mathrsfs}
\usepackage{amsfonts}
\usepackage{amsmath}
\usepackage{amsfonts,amssymb,color}
\usepackage{dsfont}
\usepackage{curves}
\usepackage{mathrsfs}
\usepackage{pifont}
\usepackage{amssymb}
\allowdisplaybreaks

\numberwithin{equation}{section}

\date{}

\def\BigRoman{\uppercase\expandafter{\romannumeral\number\count 255 }}
\def\Romannumeral{\afterassignment\BigRoman\count255=}

\begin{document}
\title{Perfect matchings and $A_{\alpha}$-spectral radius in 1-binding graphs}
\author{\small Sizhong Zhou$^{1}$\footnote{Corresponding author. E-mail address: zsz\_cumt@163.com (S. Zhou)}, Hongxia Liu$^{2}$\\
\small $1$. School of Science, Jiangsu University of Science and Technology,\\
\small Zhenjiang, Jiangsu 212100, China\\
\small $2$. School of Mathematics and Information Science, Yantai University,\\
\small Yantai, Shandong 264005, China\\
}

\maketitle
\begin{abstract}
\noindent Let $G$ be a graph with vertex set $V(G)$ and edge set $E(G)$. For $\alpha\in[0,1]$, we use $A_{\alpha}(G)$ and $\rho_{\alpha}(G)$ to denote the $A_{\alpha}$-matrix and the
$A_{\alpha}$-spectral radius of $G$, respectively. The binding number $\mbox{bind}(G)$ of $G$ is defined by $\mbox{bind}(G)=\min\left\{\frac{|N_G(X)|}{|X|}:\emptyset\neq X\subseteq V(G),N_G(X)\neq V(G)\right\}$.
If $\mbox{bind}(G)\geq1$, then $G$ is called 1-binding. A perfect matching in $G$ is a set of nonadjacent edges covering every vertex of $G$. Tutte proved that a graph $G$ of even order
has a perfect matching if and only if $o(G-S)\leq|S|$ holds for every $S\subseteq V(G)$ [W. Tutte, The factorization of linear graphs, J. Lond. Math. Soc. 22 (1947) 107--111]. In this
paper, we use Tutte's result to prove that a connected 1-binding graph $G$ of even order $n$ with $n\geq n(\alpha)$ has a perfect matching unless $G=K_1\vee(K_{n-5}\cup K_3\cup K_1)$ if
$\rho_{\alpha}(G)\geq\rho_{\alpha}(K_1\vee(K_{n-5}\cup K_3\cup K_1))$, where $n(\alpha)$ is defined as follows: $n(\alpha)=\max\{18,\frac{2+8\alpha}{1-2\alpha}\}$ if
$\alpha\in[0,\frac{1}{2})$, and $n(\alpha)=18$ if $\alpha=\frac{1}{2}$.
\\
\begin{flushleft}

(2020) Mathematics Subject Classification: 05C50, 05C70
\end{flushleft}
\end{abstract}

\section{Introduction}

Throughout this paper, we focus on finite and undirected graphs which have neither loops nor multiple edges. Let $G$ be a graph with vertex set $V(G)$ and edge set $E(G)$. We denote by $|V(G)|=n$
and $o(G)$ the order and the number of odd components of $G$, respectively. For $v\in V(G)$, the degree of $v$ in $G$ is denoted by $d_G(v)$ and write $\delta(G)=\min\{d_G(v):v\in V(G)\}$. For any
$S\subseteq V(G)$, $G[S]$ and $G-S$ denote the subgraphs of $G$ induced by $S$ and $V(G)\setminus S$, respectively. The complete graph with $n$ vertices is denoted by $K_n$. Given two vertex-disjoint
graphs $G_1$ and $G_2$, the union of $G_1$ and $G_2$ is denoted by $G_1\cup G_2$. The join $G_1\vee G_2$ is the graph derived from $G_1\cup G_2$ by adding all possible edges between $V(G_1)$ and
$V(G_2)$.

Woodall \cite{Wt} introduced the concept of the binding number. Given a subset $X\subseteq V(G)$, the set of neighbors of $X$ in $G$ is denoted by $N_G(X)$. The binding number $\mbox{bind}(G)$ of
$G$ is defined by
$$
\mbox{bind}(G)=\min\left\{\frac{|N_G(X)|}{|X|}:\emptyset\neq X\subseteq V(G), N_G(X)\neq V(G)\right\}.
$$
If $\mbox{bind}(G)\geq1$, then $G$ is called 1-binding.

Given a graph $G$ with $V(G)=\{v_1,v_2,\ldots,v_n\}$, let $A(G)=(a_{ij})_{n\times n}$ denote the adjacency matrix of $G$, where $a_{ij}=1$ if $v_iv_j\in E(G)$, and $a_{ij}=0$ otherwise. The signless
Laplacian matrix $Q(G)$ of $G$ is defined by $Q(G)=D(G)+A(G)$, where $D(G)=\mbox{diag}(d_G(v_1),d_G(v_2),\ldots,d_G(v_n))$ is the diagonal degree matrix of $G$. For any $\alpha\in[0,1]$, Nikiforov
\cite{Nm} defined the $A_{\alpha}$-matrix of a graph $G$ as $A_{\alpha}(G)=\alpha D(G)+(1-\alpha)A(G)$. We easily see that $A_0(G)=A(G)$, $2A_{\frac{1}{2}}(G)=Q(G)$ and $A_1(G)=D(G)$. In particular,
the largest eigenvalue $\rho(G)$ of $A(G)$ is called the adjacency spectral radius of $G$, the largest eigenvalue $q(G)$ of $Q(G)$ is called the signless Laplacian spectral radius of $G$ and the
largest eigenvalue $\rho_{\alpha}(G)$ of $A_{\alpha}(G)$ is called the $A_{\alpha}$-spectral radius of $G$. Obviously, $\rho_0(G)=\rho(G)$ and $2\rho_{\frac{1}{2}}(G)=q(G)$. We refer the reader to
\cite{PR,So,BM,LL,ZBW,ZZL,Za,ZZZL} for some properties of the $A_{\alpha}$-matrix of a graph.

A matching $M$ of a graph $G$ is a subset of $E(G)$ such that any two edges of $M$ have no common vertices. If every vertex of $G$ is incident to exactly one edge of $M$, then $M$ is called a perfect
matching (or 1-factor) of $G$. Tutte \cite{Tt} claimed a sufficient and necessary condition for a graph to contain a perfect matching. Based on Tutte's result, many scholars have studied the spectral
conditions for the existence of a perfect matching in a graph. O \cite{Os} provided an adjacency spectral radius condition for a graph to admit a perfect matching. Zhang and Lin \cite{ZL} showed a
distance spectral radius condition which guarantees the existence of a perfect matching in a graph. Zhao, Huang and Wang \cite{ZHW} established an $A_{\alpha}$-spectral radius condition for the
existence of a perfect matching in a graph. Fan and Lin \cite{FL} proved that a connected 1-binding graph $G$ of even order $n\geq12$ contains a perfect matching if $\rho(G)>\rho(K_1\vee(K_{n-5}\cup K_3\cup K_1))$.
Hu and Zhang \cite{HZ} verified an upper bound for the distance spectral radius to guarantee the existence of a perfect matching in a 1-binding graph. We refer the reader to \cite{Ws,Wu,LP,Zs,Zs1,Zr,Zs2,WZL,PZ,ZWH}
for other results on graph factors.

Motivated directly by \cite{Tt}, we provide a sufficient condition for a connected 1-binding graph to possess a perfect matching based on the $A_{\alpha}$-spectral radius.

\medskip

\noindent{\textbf{Theorem 1.1.}} Let $\alpha$ be a real number with $\alpha\in[0,\frac{1}{2}]$, and let $G$ be a connected 1-binding graph of even order $n$ with $n\geq n(\alpha)$, where
\[
n(\alpha)=\left\{
\begin{array}{ll}
\max\Big\{18,\frac{2+8\alpha}{1-2\alpha}\Big\},&\mbox{if} \ \alpha\in\left[0,\frac{1}{2}\right);\\
18,&\mbox{if} \ \alpha=\frac{1}{2}.\\
\end{array}
\right.
\]
If $G$ satisfies
$$
\rho_{\alpha}(G)\geq\rho_{\alpha}(K_1\vee(K_{n-5}\cup K_3\cup K_1)),
$$
then $G$ has a perfect matching, unless $G=K_1\vee(K_{n-5}\cup K_3\cup K_1)$.

\medskip

\section{Some preliminaries}

In this section, we list some essential lemmas which will be used to prove our main theorem.

\medskip

\noindent{\textbf{Lemma 2.1}} (Tutte \cite{Tt}). A graph $G$ contains a perfect matching if and only if
$$
o(G-S)\leq|S|
$$
for each subset $S\subseteq V(G)$.

\medskip

\noindent{\textbf{Lemma 2.2}} (Nikiforov \cite{Nm}). Let $H$ be a subgraph of a connected graph $G$. Then
$$
\rho_{\alpha}(G)\geq\rho_{\alpha}(H),
$$
with equality if and only if $G=H$.

\medskip

\noindent{\textbf{Lemma 2.3}} (Wu, Wang and Kang \cite{WWK}). Let $\alpha\in[0,1)$, and $\sum\limits_{i=1}^{q}n_i=n-s$. If $n_i\geq n_j$, then
$$
\rho_{\alpha}(K_s\vee(K_{n_1}\cup\cdots\cup K_{n_i}\cup\cdots\cup K_{n_j}\cup\cdots\cup K_{n_q}))<\rho_{\alpha}(K_s\vee(K_{n_1}\cup\cdots\cup K_{n_i+1}\cup\cdots\cup K_{n_j-1}\cup\cdots\cup K_{n_q})).
$$

\medskip

Let $M$ be a real $n\times n$ matrix, and let $\mathcal{N}=\{1,2,\ldots,n\}$. Given a partition $\pi:\mathcal{N}=\mathcal{N}_1\cup\mathcal{N}_2\cup\cdots\cup\mathcal{N}_r$, the matrix $M$ can be correspondingly
written as
\begin{align*}
M=\left(
  \begin{array}{cccc}
    M_{11} & M_{12} & \cdots & M_{1r}\\
    M_{21} & M_{22} & \cdots & M_{2r}\\
    \vdots & \vdots & \ddots & \vdots\\
    M_{r1} & M_{r2} & \cdots & M_{rr}\\
  \end{array}
\right),
\end{align*}
where $M_{ij}$ denotes the block of $M$ obtained by rows in $\mathcal{N}_i$ and columns in $\mathcal{N}_j$. The quotient matrix of $M$ with respect to $\pi$ is defined by the matrix $M_{\pi}=(m_{ij})_{r\times r}$,
where $m_{ij}$ is the average row sum of $M_{ij}$. If the row sum of each block $M_{ij}$ is a constant, then the partition $\pi$ is called equitable.

\medskip

\noindent{\textbf{Lemma 2.4}} (You, Yang, So and Xi\cite{YYSX}). Let $M$ be a real $n\times n$ matrix with an equitable partition $\pi$, and let $M_{\pi}$ be the corresponding quotient matrix. Then the eigenvalues
of $M_{\pi}$ are also eigenvalues of $M$. Furthermore, if $M$ is nonnegative and irreducible, then the largest eigenvalues of $M$ and $M_{\pi}$ are equal.

\medskip

A matrix $M$ is called a Hermitian matrix if $M=(\overline{M})^{T}$, where $(\overline{M})^{T}$ denotes the conjugate transpose of $M$. Obviously, a real symmetric matrix is a Hermitian matrix. The well-known Cauchy
Interlacing Theorem is given as follows.

\medskip

\noindent{\textbf{Lemma 2.5}} (Haemers \cite{Hi}). Let $M$ be a Hermitian matrix of order $s$, and let $N$ be a principal submatrix of $M$ with order $t$. If $\lambda_1\geq\lambda_2\geq\cdots\geq\lambda_s$ are the
eigenvalues of $M$ and $\mu_1\geq\mu_2\geq\cdots\geq\mu_t$ are the eigenvalues of $N$, then $\lambda_i\geq\mu_i\geq\lambda_{s-t+i}$ for $1\leq i\leq t$.

\section{The proof of Theorem 1.1}

\noindent{\it Proof of Theorem 1.1.} Suppose that a connected 1-binding graph $G$ contains no perfect matching. By Lemma 2.1, we conclude $o(G-S)\geq|S|+1$ for some subset $S\subseteq V(G)$. Note that $n$ is even.
Then $o(G-S)$ and $|S|$ admit the same parity, and so $o(G-S)\geq|S|+2$. We first prove $S\neq\emptyset$. In fact, if $S=\emptyset$, then $0=o(G)=o(G-S)\geq|S|+2=2$, which is impossible. Therefore, we infer
$|S|\geq1$.

Let $|S|=s$ and $o(G-S)=q$. Then $q\geq s+2$. We use $O_1,O_2,\ldots,O_q$ to denote the $q$ odd components in $G-S$, and write $|O_i|=n_i$ for $1\leq i\leq q$. Without loss of generality, we let
$n_q\geq n_{q-1}\geq\cdots\geq n_1$.

\medskip

\noindent{\bf Claim 1.} $n_{s+1}\geq3$.

\noindent{\it Proof.} Assume that $n_{s+1}=1$. Combining this with $n_{s+1}\geq n_s\geq\cdots\geq n_1\geq1$, we deduce $n_{s+1}=n_s=\cdots=n_1=1$. Let $X=V(O_{s+1}\cup O_s\cup\cdots\cup O_1)$. Then $|X|=s+1$ and
$N_G(X)\subseteq S$. Thus, we infer
$$
\frac{|N_G(X)|}{|X|}\leq\frac{|S|}{|X|}=\frac{s}{s+1}<1,
$$
which contradicts that $G$ is 1-binding. Hence, we obtain $n_{s+1}\geq3$. This completes the proof of Claim 1. \hfill $\Box$

\medskip

By Claim 1, we see $n_i\geq3$ for $s+1\leq i\leq q$ and $n_j\geq1$ for $1\leq j\leq s$. It is obvious that $G$ is a spanning subgraph of $G_1=K_s\vee(K_{n_{s+2}'}\cup K_{n_{s+1}}\cup K_{n_s}\cup\cdots\cup K_{n_1})$,
where $n_{s+2}'\geq n_{s+1}\geq n_s\geq\cdots\geq n_1$ are positive odd integers with $n_{s+1}\geq3$ and $\sum\limits_{i=1}^{s+1}n_i+n_{s+2}'=n-s$. By virtue of Lemma 2.2, we infer
\begin{align}\label{eq:3.1}
\rho_{\alpha}(G)\leq\rho_{\alpha}(G_1),
\end{align}
with equality following if and only if $G=G_1$.

Let $G_2=K_s\vee(K_{n-2s-3}\cup K_3\cup sK_1)$, where $n\geq2s+6$. In terms of Lemma 2.3, we conclude
\begin{align}\label{eq:3.2}
\rho_{\alpha}(G_1)\leq\rho_{\alpha}(G_2),
\end{align}
where the equality follows if and only if $(n_1,\ldots,n_s,n_{s+1},n_{s+2}')=(1,\ldots,1,3,n-2s-3)$. The following proof will be divided into two cases based on the value of $s$.

\noindent{\bf Case 1.} $s=1$.

Obviously, $G_2=K_1\vee(K_{n-5}\cup K_3\cup K_1)$. In terms of \eqref{eq:3.1} and \eqref{eq:3.2}, we get
$$
\rho_{\alpha}(G)\leq\rho_{\alpha}(K_1\vee(K_{n-5}\cup K_3\cup K_1)),
$$
where the equality occurs if and only if $G=K_1\vee(K_{n-5}\cup K_3\cup K_1)$. Observe that $K_1\vee(K_{n-5}\cup K_3\cup K_1)$ has no perfect matching. Thus, we obtain a contradiction.

\noindent{\bf Case 2.} $s\geq2$.

Recall that $G_2=K_s\vee(K_{n-2s-3}\cup K_3\cup sK_1)$. The quotient matrix of $A_{\alpha}(G_2)$ with respect to the partition $V(G_2)=V(K_s)\cup V(K_{n-2s-3})\cup V(K_3)\cup V(sK_1)$ can be given by
\begin{align*}
B_2=\left(
  \begin{array}{cccc}
  \alpha n-\alpha s+s-1 & (1-\alpha)(n-2s-3) & 3-3\alpha & (1-\alpha)s\\
  (1-\alpha)s & n+\alpha s-2s-4 & 0 & 0\\
  (1-\alpha)s & 0 & \alpha s+2 & 0\\
  (1-\alpha)s & 0 & 0 & \alpha s\\
  \end{array}
\right).
\end{align*}
The characteristic polynomial of $B_2$ is
\begin{align}\label{eq:3.3}
\varphi_{B_2}(x)=&x^{4}-\Big((1+\alpha)n+(2\alpha-1)s-3\Big)x^{3}\nonumber\\
&+\Big(\alpha n^{2}+((2\alpha^{2}+\alpha)s-2\alpha+1)n+(\alpha^{2}-\alpha-1)s^{2}-(5\alpha+4)s-6\Big)x^{2}\nonumber\\
&-\Big((2\alpha^{2}s+2\alpha)n^{2}-((-\alpha^{3}+2\alpha^{2}-2\alpha+1)s^{2}+(9\alpha^{2}-4\alpha+3)s+8\alpha+2)n\nonumber\\
&+(3\alpha^{2}-6\alpha+2)s^{3}+(11\alpha^{2}-24\alpha+8)s^{2}+(18\alpha^{2}-32\alpha+14)s+8\Big)x\nonumber\\
&+(\alpha^{3}s^{2}+2\alpha^{2}s)n^{2}-((2\alpha^{3}-2\alpha^{2}+\alpha)s^{3}+(7\alpha^{3}-\alpha^{2}-\alpha+2)s^{2}+(8\alpha^{2}+2\alpha)s)n\nonumber\\
&+(3\alpha^{3}-5\alpha^{2}+2\alpha)s^{4}+(13\alpha^{3}-16\alpha^{2}+4)s^{3}+(18\alpha^{3}-18\alpha^{2}-2\alpha+8)s^{2}+8\alpha s.
\end{align}
By virtue of Lemma 2.4 and the equitable partition $V(G_2)=V(K_s)\cup V(K_{n-2s-3})\cup V(K_3)\cup V(sK_1)$, $\rho_{\alpha}(G_2)$ is the largest root of $\varphi_{B_2}(x)=0$. Namely, $\varphi_{B_2}(\rho_{\alpha}(G_2))=0$.
Let $\theta_1=\rho_{\alpha}(G_2)\geq\theta_2\geq\theta_3\geq\theta_4$ be the four roots of $\varphi_{B_2}(x)=0$ and $Q=\mbox{diag}(s,n-2s-3,3,s)$. Then we easily see that
\begin{align*}
Q^{\frac{1}{2}}B_2Q^{-\frac{1}{2}}=\left(
  \begin{array}{cccc}
  \alpha n-\alpha s+s-1 & (1-\alpha)s^{\frac{1}{2}}(n-2s-3)^{\frac{1}{2}} & (1-\alpha)(3s)^{\frac{1}{2}} & (1-\alpha)s\\
  (1-\alpha)s^{\frac{1}{2}}(n-2s-3)^{\frac{1}{2}} & n+\alpha s-2s-4 & 0 & 0\\
  (1-\alpha)(3s)^{\frac{1}{2}} & 0 & \alpha s+2 & 0\\
  (1-\alpha)s & 0 & 0 & \alpha s\\
  \end{array}
\right)
\end{align*}
is symmetric, and
\begin{align*}
\left(
  \begin{array}{ccc}
    n+\alpha s-2s-4 & 0 & 0\\
    0 & \alpha s+2 & 0\\
    0 & 0 & \alpha s\\
  \end{array}
\right)
\end{align*}
is a submatrix of $Q^{\frac{1}{2}}B_2Q^{-\frac{1}{2}}$. Note that $Q^{\frac{1}{2}}B_2Q^{-\frac{1}{2}}$ and $B_2$ possess the same eigenvalues. Then the Cauchy Interlacing Theorem (see Lemma 2.5) leads to
$\theta_2\leq n+\alpha s-2s-4<n-5$.

Let $G_*=K_1\vee(K_{n-5}\cup K_3\cup K_1)$. The quotient matrix of $A_{\alpha}(G_*)$ with respect to the partition $V(G_*)=V(K_1)\cup V(K_{n-5})\cup V(K_3)\cup V(K_1)$ can be written as
\begin{align*}
B_*=\left(
  \begin{array}{cccc}
  \alpha n-\alpha & (1-\alpha)(n-5) & 3-3\alpha & 1-\alpha\\
  1-\alpha & n+\alpha-6 & 0 & 0\\
  1-\alpha & 0 & \alpha+2 & 0\\
  1-\alpha & 0 & 0 & \alpha\\
  \end{array}
\right).
\end{align*}
The characteristic polynomial of $B_*$ equals
\begin{align*}
\varphi_{B_*}(x)=&x^{4}-\Big((1+\alpha)n+2\alpha-4\Big)x^{3}+\Big(\alpha n^{2}+(2\alpha^{2}-\alpha+1)n+\alpha^{2}-6\alpha-11\Big)x^{2}\\
&-\Big((2\alpha^{2}+2\alpha)n^{2}+(\alpha^{3}-11\alpha^{2}-2\alpha-6)n+32\alpha^{2}-62\alpha+32\Big)x\\
&+(\alpha^{3}+2\alpha^{2})n^{2}-(9\alpha^{3}+5\alpha^{2}+2\alpha+2)n+34\alpha^{3}-39\alpha^{2}+8\alpha+12.
\end{align*}
In terms of Lemma 2.4 and the equitable partition $V(G_*)=V(K_1)\cup V(K_{n-5})\cup V(K_3)\cup V(K_1)$, $\rho_{\alpha}(G_*)$ is the largest root of $\varphi_{B_*}(x)=0$. Let $\rho_{\alpha}(G_*)=\theta$. Then
$\varphi_{B_*}(\theta)=0$.

Since $G_*=K_1\vee(K_{n-5}\cup K_3\cup K_1)$ contains $K_{n-4}$ as its proper subgraph, it follows from Lemmas 2.2 that
\begin{align}\label{eq:3.4}
\theta=\rho_{\alpha}(G_*)>\rho_{\alpha}(K_{n-4})=n-5>\theta_2.
\end{align}
Recall that $\varphi_{B_*}(\theta)=0$. Let $f(x)=(1-2\alpha)x^{3}+((2\alpha^{2}+\alpha)n+(\alpha^{2}-\alpha-1)s+\alpha^{2}-6\alpha-5)x^{2}-(2\alpha^{2}n^{2}+((\alpha^{3}-2\alpha^{2}+2\alpha-1)s+\alpha^{3}-11\alpha^{2}+6\alpha-4)n+(3\alpha^{2}-6\alpha+2)s^{2}
+(14\alpha^{2}-30\alpha+10)s+32\alpha^{2}-62\alpha+24)x+(\alpha^{3}s+\alpha^{3}+2\alpha^{2})n^{2}-((2\alpha^{3}-2\alpha^{2}+\alpha)s^{2}+(9\alpha^{3}-3\alpha^{2}+2)s+9\alpha^{3}+5\alpha^{2}+2\alpha+2)n
+(3\alpha^{3}-5\alpha^{2}+2\alpha)s^{3}+(16\alpha^{3}-21\alpha^{2}+2\alpha+4)s^{2}+(34\alpha^{3}-39\alpha^{2}+12)s+34\alpha^{3}-39\alpha^{2}+8\alpha+12$. Then a simple calculation implies that
\begin{align}\label{eq:3.5}
\varphi_{B_2}(\theta)=\varphi_{B_2}(\theta)-\varphi_{B_*}(\theta)=(s-1)f(\theta).
\end{align}

In fact, $f'(x)=3(1-2\alpha)x^{2}+2((2\alpha^{2}+\alpha)n+(\alpha^{2}-\alpha-1)s+\alpha^{2}-6\alpha-5)x-2\alpha^{2}n^{2}-((\alpha^{3}-2\alpha^{2}+2\alpha-1)s+\alpha^{3}-11\alpha^{2}+6\alpha-4)n-(3\alpha^{2}-6\alpha+2)s^{2}
-(14\alpha^{2}-30\alpha+10)s-32\alpha^{2}+62\alpha-24$. We are to prove that $f'(x)>0$ for $x\geq n-5$.

For $0\leq\alpha<\frac{1}{2}$, the symmetry axis of $f'(x)$ is $x=-\frac{(2\alpha^{2}+\alpha)n+(\alpha^{2}-\alpha-1)s+\alpha^{2}-6\alpha-5}{3(1-2\alpha)}$ and $f'(x)$ is increasing for
$x\geq-\frac{(2\alpha^{2}+\alpha)n+(\alpha^{2}-\alpha-1)s+\alpha^{2}-6\alpha-5}{3(1-2\alpha)}$.

\medskip

\noindent{\bf Claim 2.} $-\frac{(2\alpha^{2}+\alpha)n+(\alpha^{2}-\alpha-1)s+\alpha^{2}-6\alpha-5}{3(1-2\alpha)}\leq n-5$.

\noindent{\it Proof.} According to $0\leq\alpha<\frac{1}{2}$, $s\geq2$, $n\geq2s+6$ and $n\geq\frac{2+8\alpha}{1-2\alpha}$, we deduce
\begin{align*}
&3(1-2\alpha)(n-5)+(2\alpha^{2}+\alpha)n+(\alpha^{2}-\alpha-1)s+\alpha^{2}-6\alpha-5\\
=&(2\alpha^{2}-3\alpha+2)n+(1-2\alpha)n+(\alpha^{2}-\alpha-1)s+\alpha^{2}+24\alpha-20\\
\geq&(2\alpha^{2}-3\alpha+2)(2s+6)+(1-2\alpha)\Big(\frac{2+8\alpha}{1-2\alpha}\Big)+(\alpha^{2}-\alpha-1)s+\alpha^{2}+24\alpha-20\\
=&(5\alpha^{2}-7\alpha+3)s+13\alpha^{2}+14\alpha-6\\
\geq&2(5\alpha^{2}-7\alpha+3)+13\alpha^{2}+14\alpha-6\\
=&23\alpha^{2}\\
\geq&0,
\end{align*}
which implies that $-\frac{(2\alpha^{2}+\alpha)n+(\alpha^{2}-\alpha-1)s+\alpha^{2}-6\alpha-5}{3(1-2\alpha)}\leq n-5$. This completes the proof of Claim 2. \hfill $\Box$

\medskip

Recall that $f'(x)$ is increasing for $x\geq-\frac{(2\alpha^{2}+\alpha)n+(\alpha^{2}-\alpha-1)s+\alpha^{2}-6\alpha-5}{3(1-2\alpha)}$ when $0\leq\alpha<\frac{1}{2}$. If $s\geq4$, then it follows from
$x\geq n-5$, $n\geq2s+6$, $0\leq\alpha<\frac{1}{2}$ and Claim 2 that
\begin{align*}
f'(x)\geq&f'(n-5)\\
=&3(1-2\alpha)(n-5)^{2}+2((2\alpha^{2}+\alpha)n+(\alpha^{2}-\alpha-1)s+\alpha^{2}-6\alpha-5)(n-5)\\
&-2\alpha^{2}n^{2}-((\alpha^{3}-2\alpha^{2}+2\alpha-1)s+\alpha^{3}-11\alpha^{2}+6\alpha-4)n\\
&-(3\alpha^{2}-6\alpha+2)s^{2}-(14\alpha^{2}-30\alpha+10)s-32\alpha^{2}+62\alpha-24\\
=&(2\alpha^{2}-4\alpha+3)n^{2}+((-\alpha^{3}+4\alpha^{2}-4\alpha-1)s-\alpha^{3}-7\alpha^{2}+32\alpha-36)n\\
&-(3\alpha^{2}-6\alpha+2)s^{2}-(24\alpha^{2}-40\alpha)s-42\alpha^{2}-28\alpha+101\\
\geq&(2\alpha^{2}-4\alpha+3)(2s+6)^{2}+((-\alpha^{3}+4\alpha^{2}-4\alpha-1)s-\alpha^{3}-7\alpha^{2}+32\alpha-36)(2s+6)\\
&-(3\alpha^{2}-6\alpha+2)s^{2}-(24\alpha^{2}-40\alpha)s-42\alpha^{2}-28\alpha+101\\
=&(-2\alpha^{3}+13\alpha^{2}-18\alpha+8)s^{2}+(-8\alpha^{3}+34\alpha^{2}-16\alpha-6)s-6\alpha^{3}-12\alpha^{2}+20\alpha-7\\
\geq&16(-2\alpha^{3}+13\alpha^{2}-18\alpha+8)+4(-8\alpha^{3}+34\alpha^{2}-16\alpha-6)-6\alpha^{3}-12\alpha^{2}+20\alpha-7\\
=&-70\alpha^{3}+332\alpha^{2}-332\alpha+97\\
>&0.
\end{align*}
If $2\leq s\leq3$, then it follows from $x\geq n-5$, $n\geq18$, $0\leq\alpha<\frac{1}{2}$ and Claim 2 that
\begin{align*}
f'(x)\geq&f'(n-5)\\
=&(2\alpha^{2}-4\alpha+3)n^{2}+((-\alpha^{3}+4\alpha^{2}-4\alpha-1)s-\alpha^{3}-7\alpha^{2}+32\alpha-36)n\\
&-(3\alpha^{2}-6\alpha+2)s^{2}-(24\alpha^{2}-40\alpha)s-42\alpha^{2}-28\alpha+101\\
\geq&324(2\alpha^{2}-4\alpha+3)+18((-\alpha^{3}+4\alpha^{2}-4\alpha-1)s-\alpha^{3}-7\alpha^{2}+32\alpha-36)\\
&-(3\alpha^{2}-6\alpha+2)s^{2}-(24\alpha^{2}-40\alpha)s-42\alpha^{2}-28\alpha+101\\
=&(-3\alpha^{2}+6\alpha-2)s^{2}+(-18\alpha^{3}+48\alpha^{2}-32\alpha-18)s-18\alpha^{3}+480\alpha^{2}-748\alpha+425\\
=&\left\{
\begin{array}{ll}
-54\alpha^{3}+564\alpha^{2}-788\alpha+381,&\mbox{if} \ s=2\\
-72\alpha^{3}+597\alpha^{2}-790\alpha+353,&\mbox{if} \ s=3\\
\end{array}
\right.\\
>&0.
\end{align*}

For $\alpha=\frac{1}{2}$, we have $f'(x)=\frac{1}{8}((16n-20s-124)x-4n^{2}+(3s+29)n+2s^{2}+12s-8)$. If $s\geq4$, the it follows from $x\geq n-5$ and $n\geq2s+6$ that
\begin{align*}
f'(x)\geq&f'(n-5)\\
=&\frac{1}{8}((16n-20s-124)(n-5)-4n^{2}+(3s+29)n+2s^{2}+12s-8)\\
=&\frac{1}{8}(12n^{2}-(17s+175)n+2s^{2}+112s+612)\\
\geq&\frac{1}{8}(12(2s+6)^{2}-(17s+175)(2s+6)+2s^{2}+112s+612)\\
=&\frac{1}{4}(8s^{2}-26s-3)\\
>&0.
\end{align*}
If $2\leq s\leq 3$, then it follows from $x\geq n-5$ and $n\geq18$ that
\begin{align*}
f'(x)\geq&f'(n-5)\\
=&\frac{1}{8}(12n^{2}-(17s+175)n+2s^{2}+112s+612)\\
\geq&\frac{1}{8}(12\times 18^{2}-18(17s+175)+2s^{2}+112s+612)\\
=&\frac{1}{4}(s^{2}-97s+675)\\
>&0.
\end{align*}
From the above discussion, we see that $f'(x)>0$ for $x\geq n-5$. This implies that $f(x)$ is strictly increasing with respect to $x\geq n-5$. Combining this with \eqref{eq:3.4}, we conclude
\begin{align}\label{eq:3.6}
f(\theta)>&f(n-5)\nonumber\\
=&(1-2\alpha)(n-5)^{3}+\Big((2\alpha^{2}+\alpha)n+(\alpha^{2}-\alpha-1)s+\alpha^{2}-6\alpha-5\Big)(n-5)^{2}\nonumber\\
&-\Big(2\alpha^{2}n^{2}+((\alpha^{3}-2\alpha^{2}+2\alpha-1)s+\alpha^{3}-11\alpha^{2}+6\alpha-4)n+(3\alpha^{2}-6\alpha+2)s^{2}\nonumber\\
&+(14\alpha^{2}-30\alpha+10)s+32\alpha^{2}-62\alpha+24\Big)(n-5)+(\alpha^{3}s+\alpha^{3}+2\alpha^{2})n^{2}\nonumber\\
&-((2\alpha^{3}-2\alpha^{2}+\alpha)s^{2}+(9\alpha^{3}-3\alpha^{2}+2)s+9\alpha^{3}+5\alpha^{2}+2\alpha+2)n\nonumber\\
&+(3\alpha^{3}-5\alpha^{2}+2\alpha)s^{3}+(16\alpha^{3}-21\alpha^{2}+2\alpha+4)s^{2}+(34\alpha^{3}-39\alpha^{2}+12)s\nonumber\\
&+34\alpha^{3}-39\alpha^{2}+8\alpha+12\nonumber\\
=&(1-\alpha)n^{3}+\Big((3\alpha^{2}-3\alpha)s+4\alpha^{2}+8\alpha-16\Big)n^{2}\nonumber\\
&-\Big((2\alpha^{3}+\alpha^{2}-5\alpha+2)s^{2}+(4\alpha^{3}+31\alpha^{2}-50\alpha+7)s+4\alpha^{3}+52\alpha^{2}-25\alpha-79\Big)n\nonumber\\
&+(3\alpha^{3}-5\alpha^{2}+2\alpha)s^{3}+(16\alpha^{3}-6\alpha^{2}-28\alpha+14)s^{2}\nonumber\\
&+(34\alpha^{3}+56\alpha^{2}-175\alpha+37)s+34\alpha^{3}+146\alpha^{2}-202\alpha-118\nonumber\\
\triangleq&g(n).
\end{align}
For $s\geq6$, it follows from $n\geq2s+6$ and $0\leq\alpha\leq\frac{1}{2}$ that
\begin{align*}
g'(n)=&3(1-\alpha)n^{2}+2\Big((3\alpha^{2}-3\alpha)s+4\alpha^{2}+8\alpha-16\Big)n\\
&-(2\alpha^{3}+\alpha^{2}-5\alpha+2)s^{2}-(4\alpha^{3}+31\alpha^{2}-50\alpha+7)s-4\alpha^{3}-52\alpha^{2}+25\alpha+79\\
\geq&3(1-\alpha)(2s+6)^{2}+2\Big((3\alpha^{2}-3\alpha)s+4\alpha^{2}+8\alpha-16\Big)(2s+6)\\
&-(2\alpha^{3}+\alpha^{2}-5\alpha+2)s^{2}-(4\alpha^{3}+31\alpha^{2}-50\alpha+7)s-4\alpha^{3}-52\alpha^{2}+25\alpha+79\\
=&(-2\alpha^{3}+11\alpha^{2}-19\alpha+10)s^{2}+(-4\alpha^{3}+21\alpha^{2}-26\alpha+1)s-4\alpha^{3}-4\alpha^{2}+13\alpha-5\\
\geq&36(-2\alpha^{3}+11\alpha^{2}-19\alpha+10)+6(-4\alpha^{3}+21\alpha^{2}-26\alpha+1)-4\alpha^{3}-4\alpha^{2}+13\alpha-5\\
=&-100\alpha^{3}+518\alpha^{2}-827\alpha+361\\
>&0,
\end{align*}
which implies that $g(n)$ is strictly increasing for $n\geq2s+6$. Thus, we obtain
\begin{align}\label{eq:3.7}
g(n)\geq&g(2s+6)\nonumber\\
=&(1-\alpha)(2s+6)^{3}+\Big((3\alpha^{2}-3\alpha)s+4\alpha^{2}+8\alpha-16\Big)(2s+6)^{2}\nonumber\\
&-\Big((2\alpha^{3}+\alpha^{2}-5\alpha+2)s^{2}+(4\alpha^{3}+31\alpha^{2}-50\alpha+7)s+4\alpha^{3}+52\alpha^{2}-25\alpha-79\Big)(2s+6)\nonumber\\
&+(3\alpha^{3}-5\alpha^{2}+2\alpha)s^{3}+(16\alpha^{3}-6\alpha^{2}-28\alpha+14)s^{2}\nonumber\\
&+(34\alpha^{3}+56\alpha^{2}-175\alpha+37)s+34\alpha^{3}+146\alpha^{2}-202\alpha-118\nonumber\\
=&(-\alpha^{3}+5\alpha^{2}-8\alpha+4)s^{3}+(-4\alpha^{3}+14\alpha^{2}-10\alpha-4)s^{2}\nonumber\\
&+(2\alpha^{3}-30\alpha^{2}+43\alpha-15)s+10\alpha^{3}-22\alpha^{2}+20\alpha-4\nonumber\\
\triangleq&h(s).
\end{align}
According to $s\geq6$ and $0\leq\alpha\leq\frac{1}{2}$, we get
\begin{align*}
h'(s)=&3(-\alpha^{3}+5\alpha^{2}-8\alpha+4)s^{2}+2(-4\alpha^{3}+14\alpha^{2}-10\alpha-4)s\\
&+2\alpha^{3}-30\alpha^{2}+43\alpha-15\\
\geq&108(-\alpha^{3}+5\alpha^{2}-8\alpha+4)+12(-4\alpha^{3}+14\alpha^{2}-10\alpha-4)\\
&+2\alpha^{3}-30\alpha^{2}+43\alpha-15\\
=&-154\alpha^{3}+678\alpha^{2}-941\alpha+369\\
>&0,
\end{align*}
which implies that $h(s)$ is strictly increasing for $s\geq6$, and so
\begin{align*}
h(s)\geq&h(6)\\
=&-338\alpha^{3}+1382\alpha^{2}-1810\alpha+626\\
>&0.
\end{align*}
Combining this with \eqref{eq:3.7}, we conclude $g(n)>0$ when $s\geq6$.

For $2\leq s\leq5$, it follows from $n\geq18$ and $0\leq\alpha\leq\frac{1}{2}$ that
\begin{align*}
g'(n)=&3(1-\alpha)n^{2}+2\Big((3\alpha^{2}-3\alpha)s+4\alpha^{2}+8\alpha-16\Big)n\\
&-(2\alpha^{3}+\alpha^{2}-5\alpha+2)s^{2}-(4\alpha^{3}+31\alpha^{2}-50\alpha+7)s-4\alpha^{3}-52\alpha^{2}+25\alpha+79\\
\geq&3(1-\alpha)\times 18^{2}+2\Big((3\alpha^{2}-3\alpha)s+4\alpha^{2}+8\alpha-16\Big)\times 18\\
&-(2\alpha^{3}+\alpha^{2}-5\alpha+2)s^{2}-(4\alpha^{3}+31\alpha^{2}-50\alpha+7)s-4\alpha^{3}-52\alpha^{2}+25\alpha+79\\
=&-(2\alpha^{3}+\alpha^{2}-5\alpha+2)s^{2}+(-4\alpha^{3}+77\alpha^{2}-58\alpha-7)s-4\alpha^{3}+92\alpha^{2}-659\alpha+475\\
=&\left\{
\begin{array}{ll}
-20\alpha^{3}+242\alpha^{2}-755\alpha+453,&\mbox{if} \ s=2\\
-34\alpha^{3}+314\alpha^{2}-788\alpha+436,&\mbox{if} \ s=3\\
-52\alpha^{3}+384\alpha^{2}-811\alpha+415,&\mbox{if} \ s=4\\
-74\alpha^{3}+452\alpha^{2}-824\alpha+390,&\mbox{if} \ s=5\\
\end{array}
\right.\\
>&0,
\end{align*}
which implies that $g(n)$ is strictly increasing for $n\geq18$, and so
\begin{align*}
g(n)\geq&g(18)\\
=&(1-\alpha)\times18^{3}+\Big((3\alpha^{2}-3\alpha)s+4\alpha^{2}+8\alpha-16\Big)\times18^{2}\\
&-\Big((2\alpha^{3}+\alpha^{2}-5\alpha+2)s^{2}+(4\alpha^{3}+31\alpha^{2}-50\alpha+7)s+4\alpha^{3}+52\alpha^{2}-25\alpha-79\Big)\times18\\
&+(3\alpha^{3}-5\alpha^{2}+2\alpha)s^{3}+(16\alpha^{3}-6\alpha^{2}-28\alpha+14)s^{2}\\
&+(34\alpha^{3}+56\alpha^{2}-175\alpha+37)s+34\alpha^{3}+146\alpha^{2}-202\alpha-118\\
=&(3\alpha^{3}-5\alpha^{2}+2\alpha)s^{3}+(-20\alpha^{3}-24\alpha^{2}+62\alpha-22)s^{2}\\
&+(-38\alpha^{3}+470\alpha^{2}-247\alpha-89)s-38\alpha^{3}+506\alpha^{2}-2992\alpha+1952\\
=&\left\{
\begin{array}{ll}
-170\alpha^{3}+1310\alpha^{2}-3222\alpha+1686,&\mbox{if} \ s=2\\
-251\alpha^{3}+1565\alpha^{2}-3121\alpha+1487,&\mbox{if} \ s=3\\
-318\alpha^{3}+1682\alpha^{2}-2860\alpha+1244,&\mbox{if} \ s=4\\
-353\alpha^{3}+1631\alpha^{2}-2427\alpha+957,&\mbox{if} \ s=5\\
\end{array}
\right.\\
>&0.
\end{align*}
From the above discussion, we see that $g(n)>0$. Together with \eqref{eq:3.5}, \eqref{eq:3.6} and $s\geq2$, we infer
$$
\varphi_{B_2}(\theta)=(s-1)f(\theta)>(s-1)g(n)>0.
$$
Recall that $\theta>\theta_2$ (see \eqref{eq:3.4}), $\varphi_{B_2}(\rho_{\alpha}(G_2))=0$ and $\rho_{\alpha}(G_*)=\theta$. Thus, we obtain
$$
\rho_{\alpha}(G_2)<\rho_{\alpha}(G_*)=\rho_{\alpha}(K_1\vee(K_{n-5}\cup K_3\cup K_1)).
$$
Combining this with \eqref{eq:3.1} and \eqref{eq:3.2}, we conclude
$$
\rho_{\alpha}(G)\leq\rho_{\alpha}(G_1)\leq\rho_{\alpha}(G_2)<\rho_{\alpha}(K_1\vee(K_{n-5}\cup K_3\cup K_1)),
$$
which contradicts $\rho_{\alpha}(G)\geq\rho_{\alpha}(K_1\vee(K_{n-5}\cup K_3\cup K_1))$. This completes the proof of Theorem 1.1. \hfill $\Box$

\section*{Declaration of competing interest}

The authors declare that they have no known competing financial interests or personal relationships that could have appeared to influence the work reported in this paper.

\section*{Data availability}

No data was used for the research described in the article.

\section*{Acknowledgments}

The authors would like to thank the anonymous referees for their valuable comments and suggestions to improve the quality and presentation of this paper. This work was supported by the Natural
Science Foundation of Jiangsu Province (Grant No. BK20241949). Project ZR2023MA078 supported by Shandong Provincial Natural Science Foundation.

\end{document}